\documentclass[eqsecnum,floats,aps,nofootinbib,preprint]{revtex4}
\usepackage{bm,amsfonts, mathtools}
\usepackage{graphicx,color, wasysym}
\usepackage{textcomp}
\usepackage{amsmath,amssymb,latexsym,epsfig}

\def\ulamek#1#2{\mbox{\normalfont$\frac{#1}{#2}$}}

\begin{document}

\makeatletter

\title{The spherical Bessel and Struve functions and operational methods}

\author{D. Babusci}
\email{danilo.babusci@lnf.infn.it}
\affiliation{INFN - Laboratori Nazionali di Frascati, via E. Fermi, 40, IT 00044 Frascati (Roma), Italy}

\author{G. Dattoli}
\email{dattoli@frascati.enea.it}
\affiliation{ENEA - Centro Ricerche Frascati, via E. Fermi, 45, IT 00044 Frascati (Roma), Italy}

\author{K. G\'{o}rska}
\email{katarzyna.gorska@ifj.edu.pl}
\affiliation{H. Niewodnicza\'{n}ski Institute of Nuclear Physics, Polish Academy of Sciences, ul.Eljasza-Radzikowskiego 152, 
PL 31342 Krak\'{o}w, Poland}

\author{K. A. Penson}
\email{penson@lptl.jussieu.fr}
\affiliation{Laboratoire de Physique Th\'eorique de la Mati\`{e}re Condens\'{e}e,\\
Universit\'e Pierre et Marie Curie, CNRS UMR 7600\\
Tour 13 - 5i\`{e}me \'et., B.C. 121, 4 pl. Jussieu, F 75252 Paris Cedex 05, France\vspace{2mm}}

\begin{abstract}

We review some aspects of the theory of spherical Bessel functions and Struve functions by means of an operational procedure essentially of umbral nature, capable of providing the  straightforward evaluation of their definite integrals and of successive derivatives. The method we propose allows indeed the formal reduction of these family of functions to elementary ones of Gaussian type. We study the problem in general terms and present a formalism capable of providing a unifying point of view including Anger and Weber functions too. The link to the multi-index Bessel functions is also briefly discussed.

\end{abstract}

\maketitle

\section{Introduction}\label{s:intro}

The spherical Bessel and Struve functions \cite{GEAndrews00, KOldham08} are widely exploited in applications, e. g. the diffraction and scattering of radiation \cite{PMMorse53}. Albeit their properties are rather well known, we apply here a new point of view, developed in a recent series of papers (see Ref. \cite{DBabusci11}), which frames them within the operational context. This method is fairly attractive and yields a significant simplification of the underlying computations.

In the first part of the paper we deal with the theory of spherical Bessel functions, while the second  will be devoted to the study of the properties of the Struve functions. The underlying thread will be discussed in the concluding section.

The spherical Bessel functions $j_n (x)$ are linked to their cylindrical counterparts by \cite{GEAndrews00, KOldham08, ISGradshteyn07}
\begin{equation}\label{jn}
j_n (x) = \sqrt{\frac{\pi}{2\,x}}\,J_{n + \frac12} (x), \qquad\qquad (n \in \mathbb{Z}).
\end{equation}
The relevant differential equation can easily be derived as follows. We start by introducing the operators \cite{GDattoli95}
\begin{equation}\nonumber
\hat{E}_\pm = \frac{\hat{N}}{x} \mp \partial_x,
\end{equation}
where $\hat{N}$ is an operator whose action on a generic indexed function is given by $(\nu \in \mathbb{R})$
\begin{equation}\nonumber
\hat{N}\,F_\nu (x) = \nu\, F_\nu (x). 
\end{equation}
By using the relation \cite{GEAndrews00}  
\begin{equation}\nonumber
J_\nu^\prime (x) = \mp\,\left[J_{\nu \pm 1} (x) - \frac{\nu}{x}\,J_\nu (x)\right] 
\end{equation}
we obtain the following identity 
\begin{equation}\label{shop}
\hat{E}_\pm\,J_\nu (x) = J_{\nu \pm 1} (x),
\end{equation}
and, therefore
\begin{equation}\nonumber
\hat{E}_+\,\hat{E}_-\,J_\nu (x) = J_\nu (x).
\end{equation}
In the case $\nu = n + 1/2$, taking into account of definition (\ref{jn}), this equation can be rewritten in differential terms as follows
\begin{equation}\nonumber
\left(\frac{n - \ulamek{1}{2}}{x} - \partial_x\right)\,\left(\frac{n + \ulamek{1}{2}}{x} + \partial_x\right)\,\left[\sqrt{\frac{2\,x}{\pi}}\,j_n (x)\right] = \sqrt{\frac{2\,x}{\pi}}\,j_n (x)
\end{equation}
and, after a few manipulations, we end up with 
\begin{equation}\nonumber
\left\{x^2\,\partial^2_x + 2\,x\,\partial_x + \left[x^2 - n\,(n+1)\right]\right\}\,j_n (x) = 0.
\end{equation}

According to the definition (\ref{jn}), the function $j_n (x)$ can be expressed by the series
\begin{equation}\nonumber
j_n (x) = \sqrt{\frac{\pi}{2\,x}}\,\sum_{k = 0}^\infty \frac{(-1)^k}{k!\, \Gamma \left(n + k + \frac32\right)}\,\left(\frac{x}2\right)^{n + 2\,k + \frac12}
\end{equation}
which can be formally written as \cite{DBabusci11}
\begin{equation}\label{jnop}
j_n (x) = \sqrt{\frac{\pi}{2\,x}}\,\left(\hat{c}\,\frac{x}2\right)^{n + \frac12}\,e^{- \hat{c}\,(x/2)^2}\,\varphi (0),
\end{equation}
with the operator $\hat{c}$ defined by identity
\begin{equation}\nonumber
\hat{c}^{\,\alpha}\,\varphi (0) = \frac1{\Gamma(1 + \alpha)}, \qquad\qquad (\alpha \in \mathbb{R}).
\end{equation}
Although we have already exploited this method for other families of Bessel functions \cite{DBabusci11}, we get a glimpse of how the formal re-handling of the spherical Bessel functions given in Eq. (\ref{jnop}) may be useful, by considering the evaluation of the integral
\begin{equation}\nonumber
b_0 = \int_{-\infty}^\infty j_0 (x)\,\mathrm{d}x.
\end{equation}
According to Eq. (\ref{jnop}), by treating $\hat{c}$ as an ordinary constant we can reduce our problem to the evaluation of a Gaussian integral, obtaining
\begin{equation}\nonumber
b_0 = \left[\frac{\sqrt{\pi\,\hat{c}}}2\,\int_{-\infty}^{\infty} e^{-\hat{c}\,(x/2)^2} \mathrm{d}x\right]\,\varphi(0) = \frac{\sqrt{\pi\,\hat{c}}}2\,\frac{2\,\sqrt{\pi}}{\sqrt{\hat{c}}}\,\varphi (0) = \pi\,\hat{c}^{\,0}\,\varphi(0) = \pi.
\end{equation}
In the forthcoming sections we will further prove the usefulness of this procedure.

\section{Derivatives and integrals of spherical Bessel functions}
The operational definition (\ref{jnop}) of $j_{n}(x)$ allows a very straightforward derivation of the relevant generating function. We get 
\begin{align}\label{jngen}
\sum_{n = 0}^\infty \frac{t^n}{n!}\,j_n (x) &= \frac{\sqrt{\pi\,\hat{c}}}2\,\exp\left\{-\frac{\hat{c}}4\,(x^2 - 2\,x\,t)\right\}\,\varphi (0) \nonumber \\ 
&= \sqrt{\frac{\pi}2}\,\frac{J_{1/2}\,\left(\sqrt{x^2 - 2\,x\,t}\right)}{\sqrt{x^2 - 2\,x\,t}} = j_0 \left(\sqrt{x^2 - 2\,x\,t}\right).
\end{align}
The formalism of the shift operators \cite{GDattoli97} can now be usefully exploited to get further properties. By recalling indeed that 
\begin{equation}\nonumber
\exp\left\{\lambda\,\left(\frac1{z}\,\partial_z\right)\right\}\,f(z) = f \left(\sqrt{z^2 + 2\,\lambda}\right), 
\end{equation}
Eq. (\ref{jngen}) can be rewritten as 
\begin{equation}\nonumber
\sum_{n = 0}^\infty \frac{t^n}{n!}\, j_n (x) = \sum_{n = 0}^{\infty} \left. \frac{(-\xi)^m}{m!}\, \left(\frac1{x}\,\partial_x\right)^m\,j_0 (x)\right|_{\xi = x\,t}
\end{equation}
from which we obtain the well-known property \cite{KOldham08}
\begin{equation}\nonumber
j_n (x) = (-x)^n\,\left(\frac1{x}\,\partial_x\right)^n\,j_0 (x).
\end{equation}

The described approach also allows us to easily obtain a closed form for the successive derivative of the function $j_0 (x)$. By taking into account  the following identity \cite{Book}
\begin{equation}\nonumber
\partial^{ n}_x\,e^{a\,x^2} = H_n (2\,a\,x, a)\,e^{a\,x^2}
\end{equation}
where $H_n (y, z)$ are the two-variable Hermite polynomials defined by
\begin{equation}\nonumber
H_n (y, z) = n!\,\sum_{k = 0}^{[n/2]} \frac{y^{n - 2\,k}\,z^k}{k!\,(n - 2\,k)!}\,,
\end{equation}
from Eq. (\ref{jnop}) we get
\begin{align}
\partial_x^n\,j_0 (x) &= (-1)^n\,\frac{\sqrt{\pi\,\hat{c}}}2\,H_n\left(\hat{c}\,\frac{x}2, -\frac{\hat{c}}4\right)\,e^{-\hat{c}\,(x/2)^2}\,\varphi (0) \nonumber \\
&= n! \sum_{k = 0}^{[n/2]} \frac{(-1)^{n + k} (2 x)^{-k}}{k!\, (n - 2\,k)!}\,j_{n - k} (x). \nonumber
\end{align}

As a final example of application of the method, let us now consider the following integral
\begin{equation}\nonumber
b_n = \int_{-\infty}^\infty j_n (x)\,\mathrm{d}x.
\end{equation}
By using the generating function method we can write
\begin{align}\label{bt}
b (t) &= \sum_{n = 0}^\infty \frac{t^n}{n!}\,b_n = \left[\frac{\sqrt{\pi\,\hat{c}}}2\,\int_{-\infty}^\infty \exp\left\{-\frac{\hat{c}}4\,(x^2 - 2\,x\,t)\right\}\,\mathrm{d}x\right]\,\varphi (0) \nonumber \\ 
&= \pi\,\sum_{k = 0}^\infty \frac1{(k!)^2}\,\left({\frac{t}2}\right)^{2\,k} 
\end{align}
from which, by equating the coefficients of the same powers in $t$, we get 
\begin{equation}\nonumber
b_{2 n} = \frac{(2\,n)!}{(n!)^2}\,\frac{\pi}{2^{2\,n}} = \frac{\sqrt{\pi}}{n!}\,\Gamma \left(n + \frac12\right), \qquad b_{2 n +1} = 0.
\end{equation}
From Eq. (\ref{bt}) we also get
\begin{equation}\nonumber
b (t) = \pi\,I_0 (t)
\end{equation}
where $I_0 (t)$ is the modified Bessel function of order 0 \cite{KOldham08}, and, thus, as a consequence of Eq. (\ref{jngen}), the following identity results
\begin{equation}\nonumber
\int_{-\infty}^\infty j_{0}\left(\sqrt{a\,x^2 + b\,x}\right) \mathrm{d}x = \frac{\pi}{\sqrt{a}}\,I_0\left(\frac{b}{2\,\sqrt{a}}\right),
\end{equation}
that can also be viewed as an integral representation of the function $I_0$.

\section{The Struve Functions}
The Struve functions $\mathbf{H}_\alpha (x)$ are defined by the series (see Eq. (57:6:1) in Ref. \cite{KOldham08})
\begin{equation}\nonumber
\mathbf{H}_\alpha (x) = \sum_{k = 0}^\infty \frac{(-1)^k}{\Gamma\left(m + \frac32\right)\,\Gamma\left(m + \alpha + \frac32\right)}\,\left(\frac{x}2\right)^{2\,m + \alpha + 1}.
\end{equation}
The derivation of the relevant differential equation can be achieved through the shift operators defined in Eq. (\ref{shop}). By using the differentiation formula 
(see Eq. (57:10:1) in Ref. \cite{KOldham08})
\begin{equation}\nonumber
\partial_x\,\mathbf{H}_\alpha (x) = \frac12\,\left[\mathbf{H}_{\alpha - 1} (x) - \mathbf{H}_{\alpha + 1} (x) + \frac{(x/2)^\alpha}{\sqrt{\pi}\,\Gamma\left(\alpha + \frac32\right)}\right]
\end{equation}
and taking into account the recursion formula for the Struve functions
\begin{equation}\nonumber
\mathbf{H}_{\alpha + 1} (x) + \mathbf{H}_{\alpha - 1} (x) = \frac{2\,\alpha}{x}\,\mathbf{H}_{\alpha} (x) + \frac{(x/2)^\alpha}{\sqrt{\pi}\,\Gamma\left(\alpha + \frac32\right)},
\end{equation}
it is easy to show that
\begin{equation}\nonumber
\hat{E}_+\,\mathbf{H}_\alpha (x) = \mathbf{H}_{\alpha + 1} (x) - \frac{(x/2)^\alpha}{\sqrt{\pi}\,\Gamma\left(\alpha + \frac32\right)}, \qquad 
\hat{E}_-\,\mathbf{H}_{\alpha}(x) = \mathbf{H}_{\alpha - 1} (x),
\end{equation}
i.e.
\begin{equation}\nonumber
\hat{E}_+\,\hat{E}_-\,\mathbf{H}_\alpha (x) = \mathbf{H}_{\alpha} (x) - \frac{(x/2)^{\alpha - 1}}{\sqrt{\pi}\,\Gamma(\alpha + \frac12)}.
\end{equation}
In differential terms the last equation can be rewritten as follows
\begin{equation}\nonumber
\left(\frac{\alpha - 1}{x} - \partial_x\right)\,\left(\frac{\alpha}{x} + \partial_x\right)\,\mathbf{H}_\alpha (x) = \mathbf{H}_\alpha (x) - 
\frac{(x/2)^{\alpha - 1}}{\sqrt{\pi}\,\Gamma(\alpha + \frac12)}, 
\end{equation}
that can be reduced to the following non-homogeneous Bessel equation 
\begin{equation}\nonumber
\hat{B}_\alpha\,\mathbf{H}_\alpha (x) = \frac{4\,\left(x/2\right)^{\alpha + 1}}{\sqrt{\pi}\,\Gamma\left(\alpha + \frac12\right)}, \quad
\hat{B}_\alpha = (x\,\partial_x)^2 + (x^2 - \alpha^2).
\end{equation}
The properties of Struve functions can be studied by means of a slight modification of the procedure we have followed so far. To widen the perspective we consider the problem starting from 
the study of the Humbert-Bessel functions \cite{GDattoli97_1}, which are two-index Bessel-like functions defined by the series 
\begin{equation}
\label{HumBes}
J_{\mu, \nu} (x) = \sum_{k = 0}^\infty \frac{(-x)^k}{k!\,\Gamma(k + \mu + 1)\,\Gamma(k + \nu + 1)}, 
\end{equation}
and whose connection with the Struve functions is realized by the following identity
\begin{equation}\nonumber
\mathbf{H}_\alpha (x) = \left(\frac{x}2\right)^{\alpha + 1}\,\int_0^\infty e^{- s}\,J_{\frac12,\,\alpha + \frac12} \left[s\,\left(\frac{x}2\right)^2\right]\,\mathrm{d}s.
\end{equation}
According to our technique, the Humbert-Bessel functions can formally written as
\begin{equation}
\label{Jmunu}
J_{\mu, \nu} (x) = \hat{c}_1^{\,\mu}\,\hat{c}_2^{\,\nu}\,e^{- \hat{c}_1\,\hat{c}_2\,x}\,\varphi_1 (0)\,\varphi_2 (0), \qquad 
\hat{c}_i^{\,\alpha}\,\varphi_i (0) = \frac{1}{\Gamma(1 + \alpha)} \quad (i = 1,2),
\end{equation}
an expression that makes a very simple task to derive, for example, the identities reported below
\begin{align}\nonumber
\int_{-\infty}^{\infty} J_{\mu, \nu} (x^2)\,\mathrm{d}x &= \frac{\sqrt{\pi}}{\Gamma\left(\mu + \frac12\right)\,\Gamma\left(\nu + \frac12\right)},  \\  
\int_0^\infty x^{\alpha - 1}\,J_{\mu, \nu} (x)\,\mathrm{d}x &= \frac{\Gamma(\alpha)}{\Gamma(\mu - \alpha + 1)\, \Gamma(\nu - \alpha + 1)}.  \nonumber
\end{align}

The operatorial expression (\ref{Jmunu}) allows one also to easily prove that (see Eq. (2.7.2.1) in Ref. \cite{APPrudnikov98})
\begin{align}\nonumber
\int_0^\infty \mathbf{H}_\alpha (x)\,\mathrm{d}x &= \left[\Gamma\left(1 + \frac{\alpha}2\right)\,\hat{c}_1^{- (\alpha+1)/2}\,\hat{c}_2^{(\alpha - 1)/2}\,
\int_0^\infty e^{-s}\,s^{-1 - \alpha/2}\,\mathrm{d}s\right]\,\varphi_1 (0)\,\varphi_2 (0) \nonumber \\
&= - \cot\left(\alpha\,\frac{\pi}2\right) \qquad\qquad \qquad (-2 < {\mathrm Re}(\alpha) < 0), \nonumber
\end{align}
which can be further handled to get the canonical form reported in Ref. \cite{DBabusci11}, where are also discussed the conditions of validity.

A more general definition of Struve-like functions can be obtained through the definition of the following auxiliary function 
\begin{equation}\nonumber
\Delta_{\alpha, \beta, \gamma} (x) = \int_0^\infty e^{- s}\,s^{\gamma - 1}\,J_{\alpha, \beta} \left(s\,\frac{x^2}4\right)\,\mathrm{d}s
\end{equation}
i.e., by using Eq. (\ref{Jmunu})
\begin{align}\nonumber
\Delta_{\alpha, \beta, \gamma} (x) &= \left[\hat{c}_1^{\,\alpha} \hat{c}_2^{\,\beta}\,\int_0^\infty s^{\gamma - 1}\,e^{-s\,(1 + \hat{c}_1\,\hat{c}_2\,x^2/4)}\,\mathrm{d}s\right]\,\varphi_1 (0)\,\varphi_2 (0) \\[0.5\baselineskip]
&= \frac{\Gamma(\gamma)}{\Gamma(1 + \alpha)\,\Gamma(1 + \beta)}\,_{1}F_2 \left(\gamma;\,1 + \alpha, 1+\beta; \,-\frac{x^2}4\right). \nonumber
\end{align}
with $_{p}F_{q}\left((a_{p});\, (b_{q}); \, z\right)$ the generalized hypergeometric function. From the definition (\ref{HumBes}) it is also easy to show that
\begin{equation}\nonumber
\sum_{m, n = - \infty}^\infty u^m\,v^n\,J_{m, n} (x) = \exp\left\{u + v - \frac{x}{u\,v}\right\}
\end{equation}
and, therefore, the generating function of $\Delta_{\alpha, \beta, \gamma} (x)$ with respect to the indices $\alpha$ and $\beta$ is given by 
\begin{align}\nonumber
\sum_{m, n = - \infty}^\infty u^m\,v^n\,\Delta_{m, n, \gamma} (x) &= e^{u + v}\,\int_0^\infty \exp\left\{-1 - \gamma - s\,\left[1 + \frac1{u\,v}\,
\left(\frac{x}2\right)^2\right]\right\} \mathrm{d}s \nonumber \\
&= e^{u + v}\,\frac{\Gamma (\gamma)}{\left[1 + \displaystyle \frac1{u\,v}\,\left(\frac{x}2\right)^2\right]^\gamma}. \nonumber
\end{align}

The multi-index Bessel of Humbert type functions are not widespreadly known, although their use could be very advantageous in applications. As an example, we consider the case of the product of two cylindrical Bessel functions, 
which can be written as (see Eq. 8.442.1 from \cite{ISGradshteyn07})
\begin{equation}\nonumber
J_\mu (x)\,J_\nu (x) = \sum_{k = 0}^\infty \frac{(-1)^k\,\Gamma(\mu + \nu + 2 k + 1)}{k!\,\Gamma (\mu + k + 1)\,\Gamma (\nu + k + 1)\,\Gamma (\mu + \nu + k + 1)}\,
\left(\frac{x}2\right)^{2 k + \mu + \nu}
\end{equation}
i.e.
\begin{equation}
\label{Besprod}
J_\mu (x)\,J_\nu (x) = \left(\frac{x}2\right)^{\mu + \nu}\,\int_0^\infty e^{- s}\,s^{\mu + \nu}\,J_{\mu, \nu, \nu + \mu} \left(s^2\,\frac{x^2}4\right) {\mathrm d}s
\end{equation}
with
\begin{equation}\nonumber
J_{\mu, \nu, \rho} (z) = \sum_{k = 0}^\infty \frac{(-z)^k}{k!\,\Gamma(\mu + k + 1)\,\Gamma(\nu + k + 1)\,\Gamma(\rho + k + 1)}. 
\end{equation}
In the context of our operatorial method, this series can be written as 
\begin{equation}\nonumber
J_{\mu, \nu, \rho} (z) = \hat{c}_1^{\,\mu}\,\hat{c}_2^{\,\nu}\,\hat{c}_3^{\,\rho}\,e^{- \hat{c}_1\,\hat{c}_2\,\hat{c}_3\,z}\,\varphi_1 (0)\,\varphi_2 (0)\,\varphi_3 (0) 
\end{equation}
that, inserted in Eq. (\ref{Besprod}), makes an easy task to prove, for example, the following identity
\begin{align}\nonumber
\label{pino}
\int_0^\infty  \left(\frac{x}2\right)^{- \mu - \nu}J_\mu (x)\,J_\nu (x) \mathrm{d}x &= \sqrt{\pi}\,\frac{\Gamma (\mu + \nu)}{\Gamma\left(\mu + \displaystyle \frac12\right)\,
\Gamma\left(\nu + \displaystyle \frac12\right)\,\Gamma\left(\mu + \nu + \dfrac12\right)}.
\end{align}

\section{Concluding remarks}
Multi-index Bessel-like functions can also be used to express the Anger $\mathbf{J}_\nu (x)$ and Weber $\mathbf{E}_\nu (x)$ functions, defined as \cite{NIST} 
\begin{equation}\nonumber
\left(
\begin{array}{c}
\mathbf{J}_{\nu}(x) \\ 
\mathbf{E}_{\nu}(x)
\end{array}
\right) = \left(
\begin{array}{cc} 
\cos\left(\displaystyle \frac{\nu\,x}2\right) & \sin\left(\displaystyle \frac{\nu\,x}2\right) \\[10pt]
\sin\left(\displaystyle \frac{\nu\,x}2\right) & -\cos\left(\displaystyle \frac{\nu\,x}2\right) 
\end{array}
\right)\, \left(
\begin{array}{c} 
S_1 (\nu, x) \\ 
S_2 (\nu, x)
\end{array}
\right),
\end{equation}
with 
\begin{align}\nonumber
S_1 (\nu, x) &= \sum_{k = 0}^\infty \frac{(- 1)^k\,(x/2)^{2 k}}{\Gamma \left(k + \dfrac{\nu}2 + 1\right)\,\Gamma \left(k - \dfrac{\nu}2 + 1\right)} \nonumber \\[10pt]
S_2 (\nu, x) &= \sum_{k = 0}^\infty \frac{(- 1)^k\,(x/2)^{2 k + 1}}{\Gamma \left(k + \dfrac{\nu}2 + \dfrac32\right)\,
\Gamma \left(k - \dfrac{\nu}2 + \dfrac32\right)}. \nonumber
\end{align}
These functions can be expressed in the operatorial form ($\delta_{i, k}$ ($i, k = 1, 2$) is the Kronecker symbol)
\begin{equation}\nonumber
S_k (\nu, x) = \hat{c}_1^{\,(\delta_{k, 2} + \nu)/2}\,\hat{c}_2^{\,(\delta_{k, 2} - \nu)/2}\,\frac{(x/2)^{\delta_{k, 2}}}{1 + \hat{c}_1\,\hat{c}_2\,(x/2)^2}\,
\varphi_1 (0)\,\varphi_2 (0)
\end{equation}
that can be used to easily prove, for example,  the following identities 
\begin{equation}\nonumber
\int_0^\infty S_1 (\nu, x)\,{\mathrm d}x = \cos\left(\frac{\nu\,\pi}2\right), \qquad \int_0^\infty \frac{S_2 (\nu, x)}{x}\,{\mathrm d}x =  \frac1{\nu}\,\sin\left(\frac{\nu\,\pi}2\right).
\end{equation}

In this paper we have shown that the use of concepts of umbral nature from the operational calculus, combined with the properties of special functions and polynomials, can be successfully applied to the the theory of Bessel-like functions. 
The method we have suggested provides some advantages in ''practical" computations and is flexible enough to open many new perspectives, which cannot all be explored in the space of a single paper.  The main drawback is the lack of 
mathematical rigor and, therefore, the necessity of checking the obtained results resorting, for example, to numerical methods. 

Some of our results (those relevant to the integrals) can be justified on the basis of the Ramanujan master theorem \cite{GHHardy40, TAmdeberhan09, TAmdeberhan11}, originally suggested by XIX century operationalists like Glaisher 
\cite{GDattoli97}. However, we believe that the strategy one should follow, when dealing with the present method, is that of using it as a procedure to first ``guess" some specific formulae and then use them as sound conjecture to be proven 
with other conventional means.

\section{Acknowledgements}

K. A. P. acknowledges support from Agence Nationale de la Recherche (Paris, France) under Program PHYSCOMB No.~ANR-08-BLAN-243-2. 


\end{document}